
\input amstex
\documentstyle{amsppt}
\magnification 1200
\NoBlackBoxes
\input epsf

\def\refBeBr   {1}
\def\refBoPo   {2}
\def\refItSh   {3}
\def\refMa     {4}
\def\refMaWP   {5}
\def\refTopol  {6}
\def\refGAFA   {7}
\def\refPHNA   {8}
\def\refOrAa   {9}
\def\refOrPu   {10}
\def\refPoCat  {11}
\def\refPoDiss {12}
\def\refPoVNT  {13}
\def\refSturm  {14}

\def\sectFXXI   {2}
\def\sectFXX    {3}
\def\sectZapret {4}
\def\sectConstr {5}

\def\thXXI {1}
\def\thXX  {2}

\def\propPwU  {1}
\def\propPwR  {2}

\def\remMenelaus {1}
\def\remRigid    {2}

\def\figXXI    {1}
\def\figPH     {2}
\def\figSnake  {3}
\def\figLp     {4}
\def\figConstr {5}
\def\figLpII   {6}
\def\figZapret {7}
\def\figPW     {8}
\def\figPWR    {9}

\def\eqTabXXI   {1}
\def\eqTabXX    {2}
\def\eqB        {3}
\def\eqMenelaus {4}

\def\RP{\Bbb{RP}}
\def\CP{\Bbb{CP}}
\def\conj{{\bold c}}

\def\fll {\text{\rm fless}}
\def\MT  {\text{\rm MT}}
\def\Lk  {\text{\rm Lk}}
\def\Sn  {\text{\rm Sn}}
\def\Sn  {\text{\rm Fig.~\figSnake}}
\def\PH  {\exists\text{\rm ps-h}}
\def\lft {\text{\rm left}}
\def\rht {\text{\rm right}}
\def\pw  {\text{\rm pw}}

\topmatter
\title
       On arrangements of plane real quartics with respect to three lines
\endtitle
\author
       S.\,Yu.~Orevkov
\endauthor

\abstract We complete the classification of mutual arrangements of a smooth 
real algebraic or real pseudoholomorphic quartic curve and three lines
under condition that each oval of the quartic intersects the union of the lines.
This classification was started in a recent preprint by Maletto.
There is one arrangement which is realizable pseudoholomorphically but not
algebraically. It can be constructed in different ways, in particular, by a combinatorial
patchworking on an irregular triangulation. This is the first example of a combinatorial
patchworking which produces a PL curve in $\Bbb{RP}^2$ whose arrangement relative to
the coordinate axes is algebraically unrealizable.
\endabstract

\address IMT, Toulouse University, Toulouse, France \endaddress
\address Steklov Math. Inst., Moscow, Russia \endaddress

\endtopmatter

\rightheadtext{On arrangements of a plane real quartic and three lines}

\document


\head\S1. Introduction \endhead
A systematic study of mutual arrangements of transversally intersecting non-singular
real algebraic curves on $\RP^2$ up to isotopy (an analog of Hilbert's 16th
problem for reducible curves) was started by Polotovskiy in [\refPoCat], 
[\refPoDiss], [\refPoVNT] and then continued by other authors; see [\refBoPo], [\refOrAa], [\refOrPu]
and references in [\refBoPo].

When the degree increases, in contrary to the case of non-singular curves, the number of
{\it a priori} possible (and even of realizable) arrangements of reducible curves
grows much faster than the difficulty of deciding the realizability of any
particular arrangement. Therefore this problem for the degree greater than 5
was usually considered under the additional condition that all the crossing
points are real (the {\it maximal intersection condition}) and under other additional
conditions for higher degrees. Another difference
between the reducible and irreducible cases 
is that a visual identification of
arrangements for reducible curves is more difficult.

An attempt of classification of reducible curves
without the maximal intersection condition is undertaken 
in a recent preprint [\refMa].
Both realizability and non-realizability
results (as well as the identification of arrangements constructed in different ways) were
obtained using specially written computer programs. In particular, it is proven in [\refMa] that:

\roster
\item"$\bullet$" $234\le N_{411}\le 273$, where $N_{411}$ is the number of 
arrangements $(C_4,L_x,L_z)$ of a smooth real algebraic quartic curve and two lines in $\RP^2$;
                 
\item"$\bullet$" $1834\le N_{4111}^{\fll}\le 1883$,  where $N_{4111}^{\fll}$ is the number of floatless
arrangements $(C_4,L_x,L_y,L_z)$ of a smooth real algebraic quartic curve and three lines in $\RP^2$
({\it floatless}
\footnote{Ovals of $C_4$ disjoint from $L_x\cup L_y\cup L_z$ are called in [\refMa] {\it floating}.
Usually they are called {\it free}.}
means that 
$L_x\cup L_y\cup L_z$ intersects each oval of $C_4$).
\endroster

In this count we do not distinguish between arrangements obtained from each other by permuting the
lines $L_x,L_y,L_z$.
The realized arrangements are listed in [\refMaWP] in a machine-readable form.
All of them are constructed by Viro's combinatorial patchworking combined with translations
defined in [\refMa, \S6.3].
The non-realizability proofs in [\refMa] use only Harnack's inequality and Bezout theorem
for auxiliary lines, hence these results automatically extend to real pseudoholomorphic curves
(we refer to [\refGAFA, \S2], [\refPHNA, \S\S1,2] for an introduction to the subject).

The 49 floatless arrangements $(C_4,L_x,L_y,L_z)$ and 39 arrangements $(C_4,L_x,L_z)$ whose
realizability is unknown are presented in [\refMa, Figures 20 and 21]. 
In the present paper we prove that all of them except one are not realizable by subsets of real pseudoholomorphic
(in particular, real algebraic) curves. 
The remaining arrangement (as well as an arrangement obtained from it by adding a free oval; see Figure~\figPH\ below)
is realizable by real pseudoholomorphic curves but unrealizable 
by real algebraic curves. The latter fact is proven by a simplest variant of the Hilbert--Rohn--Gudkov method
(see also Remark~\remMenelaus\ at the end of \S\sectZapret).

Thus the lists in [\refMaWP] provide a complete classification of algebraically realizable
arrangements $(C_4,L_x,L_z)$ and floatless arrangements $(C_4,L_x,L_y,L_z)$.
A pseudoholomorphic analogue of these classifications is obtained by
adding one more floatless arrangement $(C_4,L_x,L_y,L_z)$.

In our opinion, the most interesting result in the present paper is an observation that
the algebraically unrealizable arrangements in Figure~\figPH\ are realizable by
a patchworking along the unique up to symmetry irregular primitive triangulation of
the triangle $[(0,0),(4,0),(0,4)]$. As far as we know, this is the first example of a
combinatorial patchworking which produces PL curves whose arrangement relative to the coordinate axes
is not realizable by a real algebraic curve of the same degree (we say that a patchworking is {\it combinatorial}
if the only input is a triangulation and a sign distribution).
Some examples of non-algebraic patchworking (not combinatorial) were found in [\refBeBr].
Notice also that an irregular patchworking produces pseudoholomorphic curves under rather general assumptions [\refItSh].

It seems plausible that the classification (started in [\refMa]) of all arrangements (not only floatless) $(C_4,L_x,L_y,L_z)$
can be considerably advanced by the methods of the present paper. However, it remains as many as 8198 open cases,
thus this should be done on a computer (maybe, by upgrading [\refMaWP]).

The author thanks I.\,V.~Itenberg and G.\,M.~Polotovskiy for useful discussions.


\head\S\sectFXXI. Classification of $(C_4,L_x,L_z)$ \endhead

Throughout the paper the notation $(C_4,L_x,L_z)$ or $(C_4,L_x,L_y,L_z)$ refers to a
real algebraic or real pseudoholomorphic (in the same almost complex structure) smooth curve of degree 4
and two or three lines transverse to each other. We say that an arrangement of one-dimensional submanifolds of $\RP^2$ is 
algebraically or pseudoholomorphically realizable
by $(C_4,L_x,L_z)$ or $(C_4,L_x,L_y,L_z)$ if it is realizable by the real loci of these algebraic or
pseudoholomorphic curves.

\proclaim{Theorem \thXXI}
The arrangements in [\refMa, Figure 21] are not pseudoholomorphically realizable by $(C_4,L_x,L_z)$.
In particular, they are not algebraically realizable by  $(C_4,L_x,L_z)$.
\endproclaim

\demo{ Proof } 
Any mutual arrangement of a nonhyperbolic (i.e. without nested ovals) quartic and a conic which have 8 common points
is either one of those shown in [\refOrPu, Figures 5--7] or is obtained from them by removing some free ovals.
In the algebraic case this fact is proved in [\refPoVNT]. A proof in the pseudoholomorphic case is given in
[\refOrPu, Proposition~1]. For each arrangement in [\refMa, Figure 21], after a perturbation of the two lines
into a conic, we obtain a non-realizable arrangement of a quartic and a conic. In the following table
we give references to Figures 5--7 in [\refOrPu] for each arrangement from [\refMa, Figure 21]:
$$
\boxed{\matrix
        5.10 &&  5.1  &&  7.2  &&  5.14 &&  5.14 &&  5.13 &&  5.5  \\
        5.4  &&  5.9  &&  5.5  &&  7.1  &&  5.12 &&  5.9  &&  5.14 \\
        5.14 &&\underline{
                 5.14}&&  5.5  &&  5.9  &&  5.13 &&  7.1  &&  5.10 \\
        7.1  &&  5.4  &&  5.2  &&  7.1  &&  5.5  &&  7.1  &&  7.1  \\
        5.4  &&  5.4  &&  5.2  &&  5.14 &&  5.4  &&  5.10 &&  5.9  \\
        5.10 &&  5.9  &&  5.4  &&  5.2  &&       &&       &&
\endmatrix}
                      \eqno(\eqTabXXI)
$$
For example, ``5.14" at the 2nd position of the 3rd row of this table means that we exclude 
the 2nd arrangement in the 3rd row in [\refMa, Figure 21] (it is also shown in Figure~\figXXI\ on the left)
as follows.
By perturbing the two lines into a conic we obtain the arrangement in Figure~\figXXI\ in the middle. It contradicts
the fact that any realizable arrangement with such non-free ovals is obtained by removing some free ovals
from the 14th item in [\refOrPu, Figure~5], which is also reproduced in Figure~\figXXI\ on the right.
\enddemo

\midinsert
\centerline{\epsfxsize 80mm\epsfbox{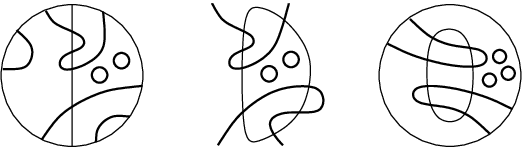}}
\botcaption{Figure \figXXI}
     Proof for the 2nd arrangement in the 3rd row in [\refMa, Fig.~21].
\endcaption
\endinsert


\head\S\sectFXX. Classification of pseudoholomorphic floatless arrangements $(C_4,L_x,L_y,L_z)$ \endhead

In this section we prove Part (b) of the following theorem.

\proclaim{ Theorem \thXX }
(a). The arrangements in [\refMa, Figure 20] are not algebraically realizable by the floatless part of $(C_4,L_x,L_y,L_z)$.
\smallskip
(b). The arrangements $(C_4,L_x,L_y,L_z)$ in Figure~\figPH\ are pseudoholomorphically realizable.
These are the only pseudoholomorphically realizable arrangements of $(C_4,L_x,L_y,L_z)$
whose floatless part is as in [\refMa, Figure 20]. Note that
the left arrangement in Figure~\figPH\ is the 44th arrangement in [\refMa, Figure 20].
\endproclaim

\midinsert
\centerline{\epsfxsize 55mm\epsfbox{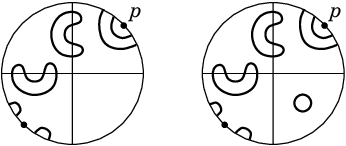}}
\botcaption{Figure \figPH}
     Algebraically non-realizable pseudoholomorphic curves
\endcaption
\endinsert

\demo{ Proof of Theorem \thXX(b) }
Notice first that if an arrangement is pseudoholomorphically realizable by $(C_4,L_x,L_y,L_z)$, then
so is its floatless part. Therefore it is enough to prove the pseudoholomorphic non-realizability
for the floatless arrangements only.

The following table as well as  [\refMa, Figure~20] has 7 rows and 7 columns.
Its entries refer to the proof for the corresponding arrangements.
$$
\boxed{
\matrix
         \Lk(z) & 7.1(y)  & \MT(z) & 7.1(z)  & 5.5(z) & \Lk(y) & \Lk(z) \\
         7.1(z) & \Sn     & 7.3(z) &  \Sn    & \MT(z) &  \Sn   & 7.1(x) \\
         7.1(y) & \Sn     & 7.1(y) & 7.1(y)  &  \Sn   &  \Sn   & 5.2(x) \\
         7.2(z) & \Sn     & 7.1(x) & \MT(z)  & 7.1(y) & 7.1(y) &  \Sn   \\
         7.2(z) & 5.5(z)  & \MT(z) & 7.1(y)  &  \Sn   & \MT(y) & \Lk(y) \\
         7.1(z) & \MT(z)  & \Lk(z) & \MT(z)  & \Lk(x) & \Lk(y) &  \Sn   \\
         7.1(x) & \PH(z)  & 7.1(y) & \MT(z)  & 7.1(x) & 7.3(z) & \MT(z)
\endmatrix}
                   \eqno(\eqTabXX)
$$

{\sl Step 1. Removal of a line; auxiliary conics.}
The entries $7.1(y)$, $7.1(z),\dots$ refer to arguments as in the proof of
Theorem~\thXXI\ applied after removal of one of the three lines. The letters
in the parentheses indicate which line should be removed: $L_x$, $L_y$, or $L_z$
(these lines are represented in [\refMa] by vertical segments, horizontal segments, and
boundary circles respectively).

An entry ``Fig.~\figSnake" refers to an unrealizability proof similar to [\refOrPu, \S2.1], which uses an auxiliary conic.
Namely, the conic passing through the five points shown in Figure~\figSnake\ has at least two intersection points
with each oval, and hence exactly two by Bezout theorem. We see in Figure~\figSnake\ that it is impossible to
trace the conic so that it crosses each line at two points.

\midinsert
\centerline{\epsfxsize 100mm\epsfbox{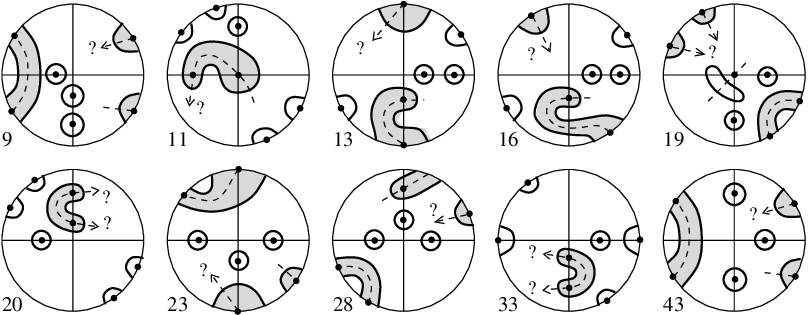}}
\botcaption{Figure \figSnake}
     Auxiliary conics.
\endcaption
\endinsert

{\sl Step 2. Braids.}
We study the other arrangements using braids. We refer to [\refTopol] and [\refGAFA, \S2] for a detailed
description of this method.
For each arrangement corresponding to the entries MT, Lk, and $\exists$ps-h 
we choose a line $L\in\{L_x, L_y,L_z\}$ (the letter
$x$, $y$, or $z$ in the table) so that there is an oval $O$ which crosses $L$ at four points.
Let $\{L_1,L_2\}$ be $\{L_x, L_y,L_z\}\setminus\{L\}$.
We choose a point $p$ on $L$ which is placed with respect to $O$ as in Figure~\figPH.
Let $\Cal L_p$ be the pencil of lines through $p$.
Then we consider the arrangements of $C_4\cup L_1\cup L_2$ on the affine plane $\RP^2\setminus L$
with respect to $\Cal L_p$ 
(called $\Cal L_p$-schemes in [\refTopol], [\refGAFA] or fiberwise arrangements in [\refOrAa]).
Up to reductions and equivalences from [\refGAFA, Corollary 2.3], there are 
three $\Cal L_p$-schemes for the arrangements nos.~1, 7, 34, two for 6, 25, 44,
and one for the other arrangements marked by MT or Lk in (2)
(the numbering corresponds to [\refMa, Table 21]). For each $\Cal L_p$-scheme we compute
a braid with 6 strands as explained in [\refGAFA, Algorithm 2.1]. An arrangement is pseudoholomorphically realizable
if and only if all the corresponding braids
are quasipositive.

All $\Cal L_p$-schemes for the arrangements marked by Lk are excluded using the linking numbers
as in [\refGAFA, \S4.5].
Those for MT are excluded by
Murasugi-Tristram inequality for braids [\refGAFA, Corollary 3.2] with the usual signature, i.e.~with $\zeta=-1$
in the notation of [\refGAFA, \S3.1].

{\sl Step 3. Realizability of Figure~\figPH.}
The arrangement no.~44 is shown in Figure~\figPH\ on the left (the corresponding entry in (\eqTabXX) is \hbox{$\exists$ps-h)}.
It is pseudoholomorphically realizable. We consider this
case in more detail (which also illustrates some underlying computations for the arrangements marked by MT and Lk).
By [\refGAFA, Corollary 2.3] it is enough to consider the two $\Cal L_p$-schemes shown in Figure~\figLp,
where $\Cal L_p$ is supposed to be the pencil of vertical lines.

\midinsert
\centerline{\epsfxsize 95mm\epsfbox{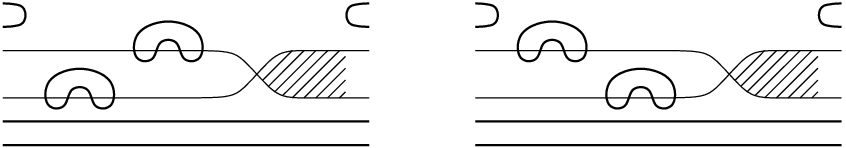}}
\botcaption{Figure \figLp}
     $\Cal L_p$-schemes for Figure \figPH\ (no.~44 in [\refMa, Fig.~20]).
\endcaption
\endinsert

\noindent
The corresponding 6-braids are, respectively,
$$
    b_1 = \bar\sigma_5 \bar\sigma_4       \sigma_5 \bar\sigma_3^4 \bar\sigma_4
          \bar\sigma_5 \bar\sigma_4^3 \bar\sigma_5 \bar\sigma_3 \Delta \qquad\text{and}\qquad
    b_2 = \bar\sigma_5 \bar\sigma_4^3 \bar\sigma_5 \bar\sigma_4 \bar\sigma_3^4
          \bar\sigma_4 \bar\sigma_5       \sigma_4 \bar\sigma_3 \Delta,                \eqno(\eqB)
$$
where $\bar\sigma_i=\sigma_i^{-1}$ and $\Delta$ is the Garside's half-twist. Both braids are
quasipositive:
$$
   b_1= \sigma_4^{\sigma_3^2\sigma_4\sigma_5\sigma_4}\sigma_2^{\sigma_3\sigma_4\sigma_5^2}
   \sigma_1^{\sigma_2\sigma_3\sigma_4\sigma_5}
   \qquad\text{and}\qquad
   b_2= \sigma_2^{\sigma_3\sigma_4\sigma_5^2}\sigma_1^{\sigma_2^2}\sigma_3^{\sigma_2\sigma_4\sigma_5},
$$
where $b^a$ denotes $a^{-1}ba$. Hence Figure \figPH(left) is pseudoholomorphically realizable.


{\sl Step 4. Uniqueness of Figure~\figPH.}
Finally, let us show that Figure~\figPH(right) is a unique pseudoholomorphically realizable arrangement
obtained by adding an oval to Figure~\figPH(left).
Proceeding as in the proof of Theorem~\thXXI, we conclude that a free oval cannot appear
anywhere else, i.e.~it might appear in the shaded zone in Figure~\figLp\ only. The corresponding
braids $b_1^+$ and $b_2^+$ are obtained by inserting
$\bar\sigma_4 \sigma_5 \bar\sigma_4\bar\sigma_5 \sigma_4$ in (\eqB) 
before $\Delta$ in $b_1$ and in $b_2$ respectively.
Then 
$b_2^+=\sigma_2^{a\sigma_5}\sigma_1^{\sigma_2a}$ ($a=\sigma_3\sigma_4\sigma_5$) is quasipositive
whereas a computation of the linking numbers shows that $b_1^+$ is not.

{\sl Step 5. A more geometric realization of Figure~\figPH.} By perturbing three lines
we may obtain a pseudoholomorphic cubic curve arranged as in Figure~\figConstr\ (cf.~[\refPHNA, before Remark~1.4]).
A further perturbation of its union with the fourth line (which is represented by the boundary circle) yields the
arrangements in Figure~\figPH.
\enddemo

\midinsert
\centerline{\epsfxsize 63mm\epsfbox{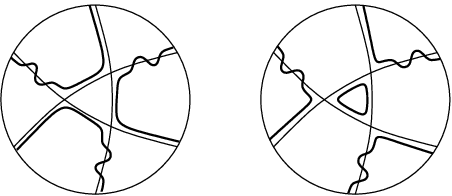}}
\botcaption{Figure \figConstr}
     Pseudoholomorphic realization of Figure \figPH.
\endcaption
\endinsert


\head\S\sectZapret. Algebraic non-realizability of the arrangements in Figure \figPH \endhead

In this section we prove Theorem \thXX(a) using a simplest version of Hilbert--Rohn--Gudkov method;
see also Remark~\remMenelaus\ below.

Suppose that there exists a real algebraic quartic curve $C$ arranged with respect to three lines
$L_x,L_y,L_z$ as in Figure~\figPH.
Then the $\Cal L_p$-scheme of its floatless part is as in Figure~\figLp.
Let us choose a line of the pencil $\Cal L_p$ tangent to an oval of $C$ as shown in Figure~\figLpII\
and roll it along the oval in the direction indicated in the picture. We do it so that the line is always
disjoint from the interior of the disk bounded by this oval (this precision is needed if
the line touches the oval at two points at some moment). During this motion the shaded domain 
shrinks and finally degenerates to a tangency point. Then we slightly push the line off from both ovals. Let
$L_0$ be the resulting line and $L=L_x\cup L_y\cup L_z\cup L_0$. Then $C\cup L$ realizes
one of the two arrangements in Figure~\figZapret, where $L_0$ is represented by the boundary circle
and the dashed free oval is or is not included.

\midinsert
\centerline{\epsfxsize 95mm\epsfbox{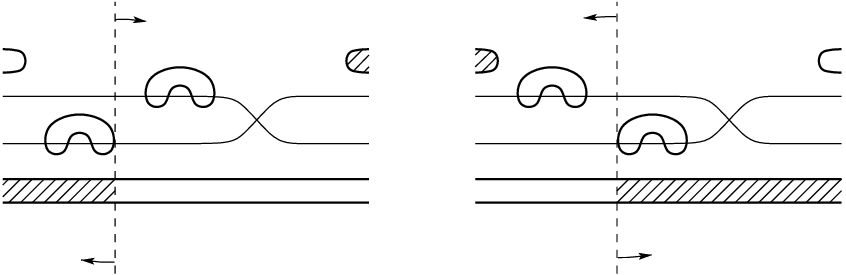}}
\botcaption{Figure \figLpII}
     Rolling a tangent line along an oval.
\endcaption
\endinsert

Note that the free oval cannot appear in the left arrangement in Figure~\figZapret\
because it cannot appear
in the left arrangement in Figure~\figLp; see Step 4 of the proof of Theorem~\thXX(b).

\midinsert
\centerline{\epsfxsize 90mm\epsfbox{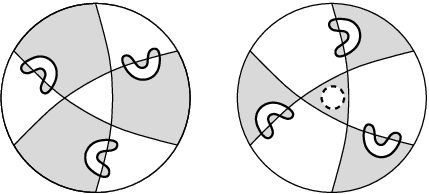}}
\botcaption{Figure \figZapret}
     Algebraic non-realizability of Figure \figPH.
\endcaption
\endinsert

Let $f(x,y,z)=0$ and $l(x,y,z)=0$ be the equations of $C_4$ and $L$ respectively.
Consider the pencil of quartics $\{C_t\}$ given by the equations $f+tl=0$.
Up to changing the sign of $l$ we may assume that the gray region in Figure~\figZapret\ grows
and the white region shrinks when $t$ varies from $t=0$ to $t>0$.
It is clear that $C_t$ cannot tend to $L$ if $C_t$ remains isotopic to $C$ for all $t>0$.
Hence there exists $t_0$ such that $C_t\cup L$ is isotopic to $C\cup L$ for $0\le t<t_0$ and
$C_{t_0}$ is singular. Let $t_1$ be such that $0<t_1-t_0\ll 1$.
It is easy to see that $C_{t_0}$ cannot be reducible: it is enough to look at the evolution of
the intersection of $C_t$ with some auxiliary lines. Hence, by the genus formula, $C_{t_0}$
has a single node.

{\it A priori} there are only three possibilities:
\roster
\item"(i)"
a node with real local branches appears in the white part of the disk bounded by a non-free oval $O$ of $C$
when two arcs collide;
\item"(ii)"
a node with imaginary local branches appears in a white region;
\item"(iii)"
the free oval of $C$ shrinks into a node with imaginary local branches.
\endroster

Case (i) is impossible because then $O$ splits into two ovals of $C_{t_1}$. This gives
the arrangement no.~29 in [\refMa, Figure~20], which is already excluded.
In Case (ii), the curve $C_{t_1}$ gets a free oval in a white region. This is impossible as we explained above.
The proof for the floatless arrangements is completed. Case (iii) is impossible because in this case
$C_{t_1}$ realizes a floatless arrangement.
Theorem~\thXX\ is proven.

\medskip\noindent
{\bf Remark \remMenelaus.} G.\,B.\,Mikhalkin observed that in order to prove the algebraic unrealizability
of the arrangements in Figure~\figZapret, instead of Hilbert--Rohn--Gudkov method, one can use
the following generalization of Menelaus' theorem. Let $C$ be a curve
of degree $n$ in $\Bbb P^2$.
Let $\xi=x/z$, $\eta=y/x$, and $\zeta=z/y$
be affine coordinates on the coordinate lines $y=0$, $z=0$, and $x=0$  respectively. 
Let $\xi_1,\dots,\xi_n$; $\eta_1,\dots,\eta_n$;
$\zeta_1,\dots,\zeta_n$ be the coordinates of the intersection points of $C$ with the coordinate lines
(here we assume that $C$ does not pass through the base points). Then
$$
     \prod_{j=1}^n \xi_j\eta_j\zeta_j = (-1)^n.              \eqno(\eqMenelaus)
$$
Indeed, let $f(x,y,z)=0$ be the equation of $C$. Then
$f(x,0,z) = a\prod(x-\xi_j z)$,
$f(0,y,z) = b\prod(y-\eta_j x)$,
$f(x,y,0) = c\prod(z-\zeta_j y)$. 
By comparing the coefficients of $x^n$, $y^n$, and $z^n$, we obtain (\eqMenelaus).
For $n=1$ this is the classical Menelaus' theorem.

If $n=4$ and $C$ realizes one of the arrangements in Figure~\figZapret, then the coordinates
can be chosen so that $L_0=\{x+y+z=0\}$ (the line represented by the boundary circle) and
all the $\xi_j$, $\eta_j$, $\zeta_j$ belong to the interval $(-1,0)$. This contradicts (\eqMenelaus).


\head\S\sectConstr. Patchworking and rigid isotopies \endhead

\midinsert
\centerline{\epsfxsize 120mm\epsfbox{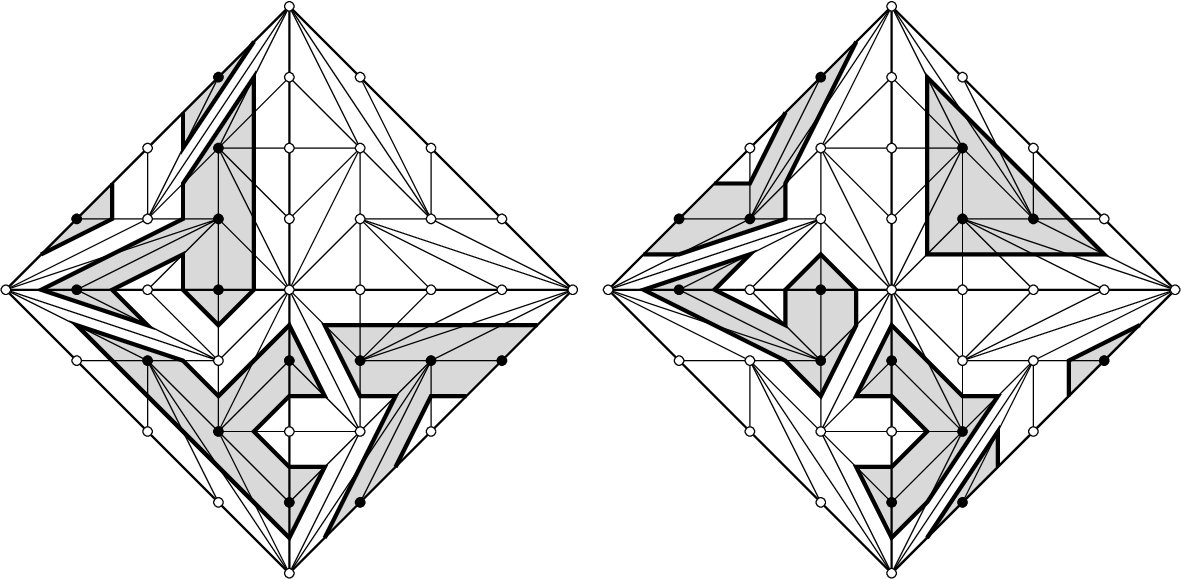}}
\botcaption{Figure \figPW}
     Patchworking of the arrangements in Figure \figPH.
\endcaption
\endinsert

The patchworking is shown in Figure \figPW\ where the signs of vertices
are represented by colors (black and white). Due to [\refItSh] it provides one more
pseudoholomorphic realization of the arrangements in Figure~\figPH.

\proclaim{ Proposition \propPwU } For each arrangement in Figure~\figPH\ there is only one (up to symmetries)
combinatorial patchworking which realizes it as a quartic curve and the three coordinate lines.
\endproclaim

\demo{ Proof } Since Figure~\figPH\ is algebraically unrealizable, the underlying triangulation must be
irregular. It is well-known that any irregular lattice triangulations of the triangle
$T=[(0,0),(4,0),(0,4)]$ is the triangulation in Figure~\figPW\ and those obtained from it by removing some
edges adjacent to points of $\partial T$ which are not vertices of $T$.
The quartic curve has four intersections with each coordinate line, hence no edge can be removed.
The quartic is disjoint from the boundary of one of the triangles bounded by the coordinate lines, hence it is enough to
consider only distributions of signs which are constant on $\partial T\cap\Bbb Z^2$. Up to symmetries, there are
only four such distributions. Two of them are shown in Figure~\figPW. It is straightforward to check that the
two others do not give an arrangement from Figure~\figPH.
\qed\enddemo

\noindent{\bf Remark \remRigid.}
It seems plausible that the patchworks in Figure \figPW\ are the only non-algebraic combinatorial
patchworks of degree 4 in $\RP^2$.
\medskip

Let $\conj:\CP^2\to\CP^2$, $(x:y:z)\mapsto(\bar x:\bar y:\bar z)$, be the involution of complex conjugation.
Let $A_0$ and $A_1$ be two real pseudoholomorphic curves in $\CP^2$, $A_k$ being $J_k$-holomorphic
for a tame $\conj$-invariant almost complex structure $J_k$.
We say that $A_0$ and $A_1$ are {\it rigidly isotopic} is there exists a $\conj$-equivariant isotopy $\{A_t\}_{t\in[0,1]}$
and a continuous family of $\conj$-anti-invariant tame almost complex structures $\{J_t\}_{t\in[0,1]}$ such that
$A_t$ is $J_t$-holomorphic for any $t$.

Let $A_{\lft}$ and $A_{\rht}$ be reducible real pseudoholomorphic curves of degree 7 realizing the floatless arrangement
in Figure~\figPH\ which correspond to the left and right $\Cal L_p$-schemes in Figure~\figLp\ respectively.
Let  $A_{\lft}^{\pw}$ and $A_{\rht}^{\pw}$ be reducible real pseudoholomorphic curves of degree 7 obtained
according to [\refItSh] from the respective patchworks in Figure~\figPW. Let also $A_{\rht}^{\pw-}$ be
the real pseudoholomorphic curve obtained from $A_{\rht}^{\pw}$ by removing the free oval (by a standard procedure
well-known in the symplectic topology).
Note that we do not know whether all these curves are defined uniquely up to rigid isotopy.

\proclaim{ Proposition \propPwR } $A_{\lft}$ and $A_{\rht}$ are not rigidly isotopic.
\endproclaim

\demo{ Proof } For any continuous deformation of the left $\Cal L_p$-scheme in Figure~\figLp\ into the right one,
there is a moment when some vertical line passes through all the three ovals, which contradicts the Bezout theorem.
\qed\enddemo

\midinsert
\centerline{\epsfxsize 68mm\epsfbox{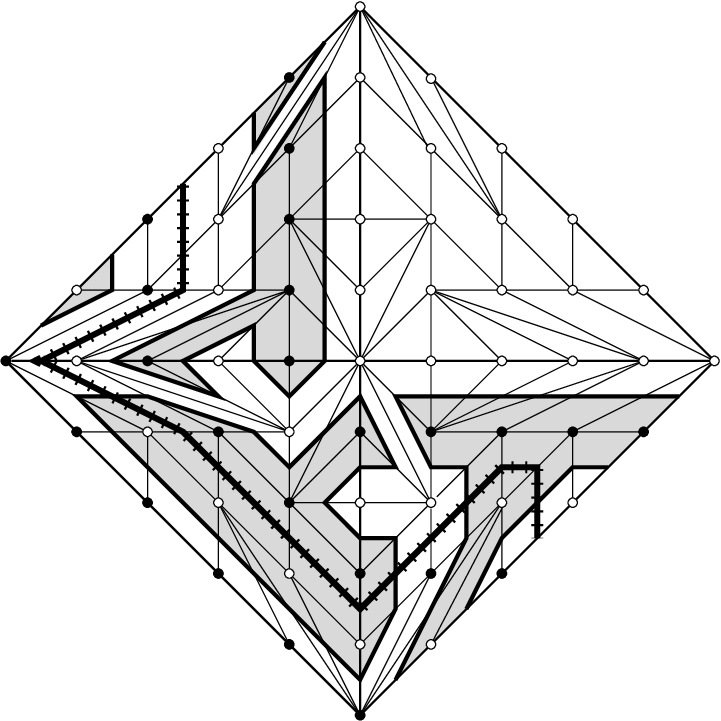}}
\vskip-1mm
\botcaption{Figure \figPWR}
     Patchworking with an additional line.
\endcaption
\endinsert

In Figure~\figPWR\ we show a patchworking of the union of a quartic curve and a line;
cf.~[\refSturm]. It gives the left arrangement in Figure~\figLp\ with a vertical line passing through
one of the two ovals (and, of course, with the line at infinity). 
This observation naturally leads to the conjecture that $A_{\lft}^{\pw}$ is not rigidly isotopic to $A_{\rht}$,
and hence, to $A_{\rht}^{\pw-}$.



\Refs

\ref\no\refBeBr\by B.~Bertrand, E.~Brugall\'e
\paper A non-algebraic pathchwork \jour Math. Z. \vol 259 \yr 2008 \issue 3 \pages 481--486 \endref

\ref\no\refBoPo\by I.\,M.~Borisov, G.\,M.~Polotovskiy 
\paper On the topology of plane real decomposable curves of degree 8
\inbook Itogi nauki i tekhniki. Sovremennye problemy matematiki. Tematicheskie obzory,
vol. 176. \yr 2020 \pages 3--18
\lang Russian \transl English transl.
\jour J. Math. Sci. \vol 275 \yr 2023 \issue 5 \pages 525--540 \endref

\ref\no\refItSh\by I.~Itenberg, E.~Shustin
\paper Combinatorial patchworking of real pseudo-holomorphic curves
\jour Turkish J. Math. \vol 26 \yr 2002 \issue 1 \pages 27--51 \endref

\ref\no\refMa\by G.~Maletto
\paper Hilbert's 16th problem for arrangements of curves on a surface
\jour\hbox to 2cm{} arxiv:2606.21449v1
\endref

\ref\no\refMaWP\by G.~Maletto \jour {\tt https://doi.org/10.5281/zenodo.20742383}, 2026\endref

\ref\no\refTopol\by S.\,Yu.~Orevkov
\paper Link theory and oval arrangements of real algebraic curves
\jour Topology \vol 38 \yr 1999 \pages 779--810 \endref

\ref\no\refGAFA\by S.\,Yu.~Orevkov
\paper Classification of flexible $M$-curves of degree 8 up to isotopy
\jour Geom. and Funct. Analysis (GAFA) \vol 12 \yr 2002 \pages 723--755 \endref

\ref\no\refPHNA\by S.\,Yu.~Orevkov
\paper Algebraically unrealizable complex orientations of plane real pseudoholomorphic curves
\jour Geom. and Funct. Analysis (GAFA) \vol 31 \yr 2021 \pages 930--947 \endref

\ref\no\refOrAa\by S.\,Yu.~Orevkov
\paper Arrangements of a plane M-sextic with respect to a line
\jour Algebra i analiz \vol 34 \yr 2022 \issue 1 \pages 123--143
\lang Russian \transl English transl.
\jour St. Petersburg Math. J. \vol 34 \yr 2023 \pages 93--107 \endref 

\ref\no\refOrPu\by S.\,Yu.~Orevkov, N.\,D.~Puchkova
\paper On mutual arrangements of a plane real curve relative to
       an $M$-quartic with an oval-snake
\jour  Mat. Zametki \vol 119 \yr 2026 \issue 2 \pages 253--265 \lang Russian
\transl English transl.
\jour  Math. Notes \vol 119 \yr 2026 \issue 2 \pages 289--298 \endref

\ref\no\refPoCat\by G.\,M.~Polotovski\u\i
\paper A catalogue of $M$-decomposing curves of sixth order
\jour Soviet Math. Doklady \vol 18 \yr 1977 \pages 1242--1246 \endref

\ref\no\refPoDiss\by G.\,M.~Polotovskiy
\paper Topological classification of decomposing curves of the 6th order
\jour  Diss. Ph.D., Gorky State Univ., Gorky, 1979 (Russian) \endref

\ref\no\refPoVNT\by G.\,M.~Polotovskiy
\paper Complete classification of $M$-decomposing curves of 6 order in the real projective plane
	\jour Dep. VINITI 20.04.1978. DEP. №1349-78, 103 p (Russian) \endref

\ref\no\refSturm\by B.~Sturmfels
\paper Viro's theorem for complete intersections
\jour Ann. Scuola Norm. Sup. Pisa Cl. Sci. (4) \vol 21 \yr 1994 \pages 377--386\endref

\endRefs
\enddocument